\documentclass{article}

\usepackage[latin1]{inputenc}
\usepackage{amsfonts,amssymb,upgreek}

\usepackage{indentfirst}

\newtheorem{teo}{Theorem}[section]
\newtheorem{lema}[teo]{Lemma}
\newtheorem{cor}[teo]{Corolary}
\newtheorem{df}[teo]{Definition}
\newtheorem{obs}[teo]{Remark}
\newtheorem{ex}[teo]{Example}
\newcommand{\bbr}{\mathbb{R}}

\begin{document}

\title{Stochastic symmetries and transformations of stochastic differential equations}
\author{ Pedro J. Catuogno\thanks{Departamento de Matemática, Universidade Estadual de Campinas, Campinas, SP, Brazil}\,\, and Luis R. Lucinger\thanks{luis@ime.unicamp.br}}
\date{\today}

\maketitle

\begin{abstract}
In this article, we introduce the notion of stochastic symmetry of a differential equation. It consists in a stochastic flow that acts over a solution of a differential equation and produces another solution of the same equation. In the ordinary case, we give necessary conditions in order to obtain such symmetries. These conditions involve the infinitesimal generator of the flow and the coefficients of the equation.

Moreover, we show how to obtain necessary conditions in order to find an application that transforms a stochastic differential equation that one would like to solve into a target equation that one previously know how to solve.
\end{abstract}

\section{Introduction}
Lie's theory of symmetries of differential equations have been studied for a long time. It is a very powerful tool if one wants to solve a differential equation explicitly. Although there is a general theory for deterministic equations (\cite{stefani}, \cite{olver}, \cite{bluman1}), the study of stochastic differential equations (SDEs) using symmetries has began only recently (\cite{albeverio}, \cite{gaeta}, \cite{wafo}, \cite{unal}, \cite{meleshko}). A major step in this theory is to obtain the determining equations. They provide conditions in order to determine if a group action is a symmetry for a differential equation.

In this article, we obtain determining equations in the case of stochastic symmetries (see Theorem \ref{teo1}). These conditions are suitable for deterministic and for stochastic equations. Further, we give conditions in order to find transformations between SDEs.

In section 2, we first introduce the theory of standard symmetries for the SDEs. Then we define what is a stochastic symmetry, obtain the determining equations for stochastic symmetries of ordinary SDEs, and give some examples.

In section 3, we deal with transformations of SDEs. It consists of finding an application $\mu$ that transforms a given SDE (that one wants to find a solution) into a target SDE (that one already knows the solution). Doing so, one can obtain the solution of a desired SDE by inverting $\mu$. This theory was developed by G. W. Bluman, A. F. Cheviakov and S. C. Anco (see \cite{bluman}, chapter 2) for deterministic differential equations. As a simple application, Bluman's theory can be used to find, systematically, the so-famous Hopf-Cole transformation that relates the heat equation and Burgers' equation. In this section, we show how to obtain $\mu$ in the case of SDEs and ilustrate the method with an example.

\section{Stochastic symmetries}

Consider the following SDE on $\bbr^n$ in the sense of Itô: \begin{equation}\label{ede}dX(t)=f(t,X(t))dt+g(t,X(t))dW(t),\end{equation} where $f$ takes value in $\bbr^n$, $W(t)$ is a Wiener process in $\bbr^m$ and $g(t,x):\bbr^m\to\bbr^n$, for each $(t,x)\in\bbr_+\times\bbr^n$. Recently, some authors have extended Lie's theory of symmetries for such equations. Let us briefly introduce that.

Let $\textbf{v}=\uptau(t)\frac{\partial}{\partial t}+\phi(t,x)\frac{\partial}{\partial x}$ be a vector field in the $(t,x)$-space. Such vector field is a symmetry for a differential equation if its flow leaves the solutions invariant. This means that for every solution of (\ref{ede}), its perturbation by the flow of $\textbf{v}$ is still a solution of (\ref{ede}). Using the coordinates of the flow,
\begin{eqnarray}\label{fluxo11}\beta_\varepsilon(t)&=&t+\int_0^\varepsilon\uptau(\beta_r(t),F_r(t,x))dr\\\label{fluxo12}F_\varepsilon(t,x)&=&x+\int_0^\varepsilon\phi(\beta_r(t),F_r(t,x))dr,\end{eqnarray}
one can find the so-called determining equations. They are: \begin{eqnarray}\label{condicaosimetria1}f_t\uptau+\uptau_tf+f_x\phi&=&\phi_t+\phi_x f+\frac{1}{2}\phi_{xx} g^2\\\label{condicaosimetria2}g_t\uptau+\frac{1}{2}\uptau_tg+g_x\phi&=&\phi_xg.\end{eqnarray}

Conditions (\ref{condicaosimetria1}) and (\ref{condicaosimetria2}) are very important. They provide necessary conditions for a flow to be a symmetry of a differential equation (or they can be used to define the admitted Lie group of transformations). Note that, even for stochastic equations, the flow is considered to be deterministic. So, it is natural to ask what happens if we change it to a stochastic flow. In what follows, we are going to consider that the perturbation of the flow in the $x$-variable (spatial variable) is stochastic. More precisely, we preserve equation (\ref{fluxo11}), but change (\ref{fluxo12}):
\begin{eqnarray}\label{fluxo21}\overline{t}=\beta_\varepsilon(t)&=&t+\int_0^\varepsilon\uptau(\beta_r(t),F_r(t,x))dr\\\label{fluxo22}\overline{x}=F_\varepsilon(t,x)&=&x+\int_0^\varepsilon\phi(\beta_r(t),F_r(t,x))dr\\\nonumber&&+\int_0^\varepsilon\tilde{\phi}(\beta_r(t),F_r(t,x))dW(r).\end{eqnarray} Here we are considering that the vector field $\textbf{v}$ has the form \begin{equation}\textbf{v}=\textbf{v}^D+\textbf{v}^S,\end{equation} where $\textbf{v}^D=\uptau(t,x)\frac{\partial}{\partial t}+\phi(t,x)\frac{\partial}{\partial x}$ is its deterministic part and $\textbf{v}^S=\tilde{\phi}(t,x)\frac{\partial}{\partial x}$ its stochastic part. This means that if $\Psi=\Psi_\varepsilon(t,x)$ is the flow of $\textbf{v}$, then \begin{equation}\Psi_\varepsilon(t,x)=(t,x)+\int_0^\varepsilon \textbf{v}^D(\Psi_r(t,x))dr+\int_0^\varepsilon\textbf{v}^S(\Psi_r(t,x))dW(r),\end{equation} which is the same as equations (\ref{fluxo21}) and (\ref{fluxo22}). An alternative notation is 
\begin{equation}\label{notation}\textbf{v}=\left[\uptau(t,x)\frac{\partial}{\partial t}+\phi(t,x)\frac{\partial}{\partial x}\right]^D+\left[\tilde{\phi}(t,x)\frac{\partial}{\partial x}\right]^S.\end{equation}

\begin{df}
A vector field $\textbf{\emph{v}}=\textbf{\emph{v}}^D+\textbf{\emph{v}}^S$ is a stochastic symmetry of a differential equation if every solution of such equation is mapped to another solution of the same equation by the flow associated to $\textbf{\emph{v}}$.
\end{df}

What we are going to do now is to obtain the determining equations for stochastic symmetries in the ordinary case. So, we want to know when the flow of a vector field $\textbf{v}$ (with a stochastic part in the spatial variable) keeps invariant the solutions of (\ref{ede}).

Let $X=X(t)$ be a solution of (\ref{ede}). Its perturbation by the flow of $\textbf{v}$ is $\overline{X}(\overline{t})$. Hence, $\textbf{v}$ will be a symmetry for (\ref{ede}) if 
\begin{equation}\label{ede2}d\overline{X}(\overline{t})=f(\overline{t},\overline{X}(\overline{t}))d\overline{t}+g(\overline{t},\overline{X}(\overline{t}))d\overline{W}(\overline{t}),\end{equation}
 where $\overline{W}$ is the Wiener process $W$ transformed by $t\mapsto\overline{t}$. By now we have to calculate the differential of $\overline{X}(\overline{t})$ using the Itô's formula and compare the result with the right hand side of (\ref{ede2}).
Note that the time change $t\mapsto\overline{t}$ can not be arbitrary. We are going to consider the time change from \cite{oksendal} (chapter 8.5). So, 
\begin{equation}\label{tempo}\overline{t}=\beta_\varepsilon(t)=\int_0^t\eta_\varepsilon^2(s)ds,\end{equation} 
with 
\begin{equation}\alpha_\varepsilon(t)=\inf\{s\geq0:\beta_\varepsilon(s)>t\}\end{equation} 
being the inverse function of $\beta$ (in the $t$ variable). Note that since $\left.\frac{\partial\beta_\varepsilon}{\partial t}(t)\right|_{\varepsilon=0}>0$, we can find a function $\eta_\varepsilon\neq0$ satisfying (\ref{tempo}) in a neighborhood of $\varepsilon=0$. Using this, the action of the flow on the solution will be given by 
\begin{equation}X(t)\mapsto \overline{X}(\overline{t})=F_\varepsilon(\alpha_\varepsilon(\overline{t}),X(\alpha_\varepsilon(\overline{t})))=F_\varepsilon(t,X(t)).\end{equation} 
For $\varepsilon$ fixed, we use Itô's formula to obtain 
\begin{eqnarray}\nonumber F_\varepsilon(t,X(t))&=&F_\varepsilon(0,X(0))+\int_0^t\left(\frac{\partial F_\varepsilon}{\partial t}+\frac{\partial F_\varepsilon}{\partial x}f+\frac{1}{2}\frac{\partial^2F_\varepsilon}{\partial x^2}g^2\right)(s,X(s))ds\\\label{eq1}&&+\int_0^t\left(\frac{\partial F_\varepsilon}{\partial x}g\right)(s,X(s))dW(s).\end{eqnarray} 
By changing variables, the first integral is equal to 
\begin{equation}\int_0^{\overline{t}}\left(\frac{\partial F_\varepsilon}{\partial t}+\frac{\partial F_\varepsilon}{\partial x}f+\frac{1}{2}\frac{\partial^2F_\varepsilon}{\partial x^2}g^2\right)(\alpha_\varepsilon(s),X(\alpha_\varepsilon(s)))\frac{\partial\alpha_\varepsilon}{\partial t}(s)ds\end{equation} 
which is the same as 
\begin{equation}\int_0^{\overline{t}}\frac{1}{\eta_\varepsilon^2(\alpha_\varepsilon(s))}\left(\frac{\partial F_\varepsilon}{\partial t}+\frac{\partial F_\varepsilon}{\partial x}f+\frac{1}{2}\frac{\partial^2F_\varepsilon}{\partial x^2}g^2\right)(\alpha_\varepsilon(s),X(\alpha_\varepsilon(s)))ds.\end{equation}
And according to \cite{oksendal}, the second integral of (\ref{eq1}) is equal to 
\begin{equation}\int_0^{\overline{t}}\frac{1}{\eta_\varepsilon(\alpha_\varepsilon(s))}\left(\frac{\partial F_\varepsilon}{\partial x}g\right)(\alpha_\varepsilon(s),X(\alpha_\varepsilon(s)))d\overline{W}(s).\end{equation}
Hence, equation (\ref{eq1}) can be rewritten as
\begin{eqnarray}\nonumber \overline{X}(\overline{t})\!\!&\!\!=\!\!&\!\!\overline{X}(0)+\int_0^{\overline{t}}\frac{1}{\eta_\varepsilon^2(\alpha_\varepsilon(s)}\left(\frac{\partial F_\varepsilon}{\partial t}+\frac{\partial F_\varepsilon}{\partial x}f+\frac{1}{2}\frac{\partial^2F_\varepsilon}{\partial x^2}g^2\right)(\alpha_\varepsilon(s),X(\alpha_\varepsilon(s)))ds\\\label{eq2}&&+\int_0^{\overline{t}}\frac{1}{\eta_\varepsilon(\alpha_\varepsilon(s))}\left(\frac{\partial F_\varepsilon}{\partial x}g\right)(\alpha_\varepsilon(s),X(\alpha_\varepsilon(s)))d\overline{W}(s).\end{eqnarray}
Comparing equations (\ref{ede2}) and (\ref{eq2}), we see that $\textbf{v}$ is a symmetry for (\ref{ede}) if and only if for all $\varepsilon$,
\begin{eqnarray}f(\overline{t},\overline{X}(\overline{t}))&=&\frac{1}{\eta_\varepsilon^2(t)}\left(\frac{\partial F_\varepsilon}{\partial t}+\frac{\partial F_\varepsilon}{\partial x}f+\frac{1}{2}\frac{\partial^2F_\varepsilon}{\partial x^2}g^2\right)(t,X(t))\\g(\overline{t},\overline{X}(\overline{t}))&=&\frac{1}{\eta_\varepsilon(t)}\left(\frac{\partial F_\varepsilon}{\partial x}g\right)(t,X(t)),\end{eqnarray}
that is,
\begin{eqnarray}\label{eq3}\!\!f(\beta_\varepsilon(t),F_\varepsilon(t,X(t)))\eta_\varepsilon^2(t)&\!\!\!=\!\!\!&\left(\frac{\partial F_\varepsilon}{\partial t}+\frac{\partial F_\varepsilon}{\partial x}f+\frac{1}{2}\frac{\partial^2F_\varepsilon}{\partial x^2}g^2\right)\!(t,X(t))\\\label{eq4}\!\!g(\beta_\varepsilon(t),F_\varepsilon(t,X(t)))\eta_\varepsilon(t)&\!\!\!=\!\!\!&\left(\frac{\partial F_\varepsilon}{\partial x}g\right)\!(t,X(t)).\end{eqnarray}
At this ponint, we calculate the differential of equations (\ref{eq3}) and (\ref{eq4}) with respect to $\varepsilon$ (note that for standard symmetries, we would calculate the derivative). For the right hand side, we use equation (\ref{fluxo22}) and the chain rule to obtain the following Lemmas.

\begin{lema}
$$d_\varepsilon\left(\displaystyle\frac{\partial F_\varepsilon}{\partial t}+\displaystyle\frac{\partial F_\varepsilon}{\partial x}\cdot f+\displaystyle\frac{1}{2}\displaystyle\frac{\partial^2F_\varepsilon}{\partial x^2}\cdot g^2\right)(t,X(t))\qquad\qquad\qquad\qquad\qquad\qquad\qquad\qquad$$
\begin{eqnarray*}&=&\left\{\frac{\partial\phi}{\partial t}(t,F_\varepsilon(t,X(t)))+\frac{\partial\phi}{\partial x}(t,F_\varepsilon(t,X(t)))\cdot\frac{\partial F_\varepsilon}{\partial t}(t,X(t))\right.\\&&\qquad+\frac{\partial\phi}{\partial x}(t,F_\varepsilon(t,X(t)))\cdot\frac{\partial F_\varepsilon}{\partial x}(t,X(t))\cdot f(t,X(t))\\&&\qquad+\frac{1}{2}\frac{\partial^2\phi}{\partial x^2}(t,F_\varepsilon(t,X(t)))\cdot\frac{\partial F_\varepsilon}{\partial x}(t,X(t))^2\cdot g(t,X(t))^2\\&&\qquad\left.+\frac{1}{2}\frac{\partial\phi}{\partial x}(t,F_\varepsilon(t,X(t)))\cdot\frac{\partial^2 F_\varepsilon}{\partial x^2}(t,X(t))\cdot g(t,X(t))^2\right\}d\varepsilon\\&&+\left\{\frac{\partial\tilde{\phi}}{\partial t}(t,F_\varepsilon(t,X(t)))+\frac{\partial\tilde{\phi}}{\partial x}(t,F_\varepsilon(t,X(t)))\cdot\frac{\partial F_\varepsilon}{\partial t}(t,X(t))\right.\\&&\qquad+\frac{\partial\tilde{\phi}}{\partial x}(t,F_\varepsilon(t,X(t)))\cdot\frac{\partial F_\varepsilon}{\partial x}(t,X(t))\cdot f(t,X(t))\\&&\qquad+\frac{1}{2}\frac{\partial^2\tilde{\phi}}{\partial x^2}(t,F_\varepsilon(t,X(t)))\cdot\frac{\partial F_\varepsilon}{\partial x}(t,X(t))^2\cdot g(t,X(t))^2\\&&\qquad\left.+\frac{1}{2}\frac{\partial\tilde{\phi}}{\partial x}(t,F_\varepsilon(t,X(t)))\cdot\frac{\partial^2 F_\varepsilon}{\partial x^2}(t,X(t))\cdot g(t,X(t))^2\right\}dW(\varepsilon).\end{eqnarray*}
\end{lema}

\begin{lema}
$$d_\varepsilon\left(\frac{\partial F_\varepsilon}{\partial x}\cdot g\right)(t,X(t))\qquad\qquad\qquad\qquad\qquad\qquad\qquad\qquad\qquad\qquad$$
\begin{eqnarray*}&=&\left\{\frac{\partial\phi}{\partial x}(t,F_\varepsilon(t,X(t)))\cdot\frac{\partial F_\varepsilon}{\partial x}(t,X(t))\cdot g(t,X(t))\right\}d\varepsilon\\&&+\left\{\frac{\partial\tilde{\phi}}{\partial x}(t,F_\varepsilon(t,X(t)))\cdot\frac{\partial F_\varepsilon}{\partial x}(t,X(t))\cdot g(t,X(t))\right\}dW(\varepsilon).\end{eqnarray*}
\end{lema}

\vskip0.5cm
For the left hand side of equations (\ref{eq3}) and (\ref{eq4}), we have

\begin{lema}\label{le}
\begin{eqnarray*}d_\varepsilon f(\beta_\varepsilon(t),F_\varepsilon(t,X(t)))&=&\left\{\frac{\partial f}{\partial t}\uptau+\frac{\partial f}{\partial x}\cdot \phi+\frac{1}{2}\frac{\partial^2 f}{\partial x^2}\cdot\tilde{\phi}^2\right\}d\varepsilon\\&&+\frac{\partial f}{\partial x}\cdot\tilde{\phi}\,dW(\varepsilon)\end{eqnarray*}
and
\begin{eqnarray*}d_\varepsilon g(\beta_\varepsilon(t),F_\varepsilon(t,X(t)))&=&\left\{\frac{\partial g}{\partial t}\uptau+\frac{\partial g}{\partial x}\cdot \phi+\frac{1}{2}\frac{\partial^2 g}{\partial x^2}\cdot\tilde{\phi}^2\right\}d\varepsilon\\&&+\frac{\partial g}{\partial x}\cdot\tilde{\phi}\,dW(\varepsilon).\end{eqnarray*}
\end{lema}
\noindent\textit{Proof:} It follows from Itô's formula to the functions $f$ and $g$ composed with $\varepsilon\mapsto(\beta_\varepsilon(t),F_\varepsilon(t,X(t)))$, together with equations (\ref{fluxo21}) and (\ref{fluxo22}).$\hfill\square$

\begin{cor}
\begin{eqnarray*}d_\varepsilon\,\eta_\varepsilon^2(t)f(\beta_\varepsilon(t),Y_\varepsilon(t))&=&\left\{\eta_\varepsilon^2\left[\frac{\partial f}{\partial t}\frac{\partial\beta_\varepsilon}{\partial\varepsilon}+\frac{\partial f}{\partial x}\cdot \phi+\frac{1}{2}\frac{\partial^2 f}{\partial x^2}\cdot\tilde{\phi}^2\right]+f2\eta_\varepsilon\frac{\partial\eta_\varepsilon}{\partial\varepsilon}\right\}d\varepsilon\\&&+\eta_\varepsilon^2\frac{\partial f}{\partial x}\cdot\tilde{\phi}\,dW(\varepsilon)\end{eqnarray*}
and
\begin{eqnarray*}d_\varepsilon\,\eta_\varepsilon(t)g(\beta_\varepsilon(t),Y_\varepsilon(t))&=&\left\{\eta_\varepsilon\left[\frac{\partial g}{\partial t}\frac{\partial\beta_\varepsilon}{\partial\varepsilon}+\frac{\partial g}{\partial x}\cdot \phi+\frac{1}{2}\frac{\partial^2 g}{\partial x^2}\cdot\tilde{\phi}^2\right]+g\frac{\partial\eta_\varepsilon}{\partial\varepsilon}\right\}d\varepsilon\\&&+\eta_\varepsilon\frac{\partial g}{\partial x}\cdot\tilde{\phi}\,dW(\varepsilon).\end{eqnarray*}
\end{cor}
\noindent\textit{Proof:} Simple use of the product rule of stochastic calculus and Lemma \ref{le}.$\hfill\square$

\vskip 0.5cm
Hence, we obtained the differential of equations (\ref{eq3}) and (\ref{eq4}). Comparing the result, we arrive at the following four equations:
$$\eta_\varepsilon^2\left[\frac{\partial f}{\partial t}\frac{\partial\beta_\varepsilon}{\partial\varepsilon}+\frac{\partial f}{\partial x}\cdot \phi+\frac{1}{2}\frac{\partial^2 f}{\partial x^2}\cdot\tilde{\phi}^2\right]+f2\eta_\varepsilon\frac{\partial\eta_\varepsilon}{\partial\varepsilon}\qquad\qquad\qquad\qquad\qquad\qquad\qquad$$\begin{eqnarray*}=&\frac{\partial\phi}{\partial t}(t,F_\varepsilon(t,X(t)))+\frac{\partial\phi}{\partial x}(t,F_\varepsilon(t,X(t)))\cdot\frac{\partial F_\varepsilon}{\partial t}(t,X(t))\\&+\frac{\partial\phi}{\partial x}(t,F_\varepsilon(t,X(t)))\cdot\frac{\partial F_\varepsilon}{\partial x}(t,X(t))\cdot f(t,X(t))\\&+\frac{1}{2}\frac{\partial^2\phi}{\partial x^2}(t,F_\varepsilon(t,X(t)))\cdot\frac{\partial F_\varepsilon}{\partial x}(t,X(t))^2\cdot g(t,X(t))^2\\&+\frac{1}{2}\frac{\partial\phi}{\partial x}(t,F_\varepsilon(t,X(t)))\cdot\frac{\partial^2 F_\varepsilon}{\partial x^2}(t,X(t))\cdot g(t,X(t))^2,\end{eqnarray*}
\begin{eqnarray*}\eta_\varepsilon^2\frac{\partial f}{\partial x}\cdot\tilde{\phi}&=&\frac{\partial\tilde{\phi}}{\partial t}(t,F_\varepsilon(t,X(t)))+\frac{\partial\tilde{\phi}}{\partial x}(t,F_\varepsilon(t,X(t)))\cdot\frac{\partial F_\varepsilon}{\partial t}(t,X(t))\\&&+\frac{\partial\tilde{\phi}}{\partial x}(t,F_\varepsilon(t,X(t)))\cdot\frac{\partial F_\varepsilon}{\partial x}(t,X(t))\cdot f(t,X(t))\\&&+\frac{1}{2}\frac{\partial^2\tilde{\phi}}{\partial x^2}(t,F_\varepsilon(t,X(t)))\cdot\frac{\partial F_\varepsilon}{\partial x}(t,X(t))^2\cdot g(t,X(t))^2\\&&+\frac{1}{2}\frac{\partial\tilde{\phi}}{\partial x}(t,F_\varepsilon(t,X(t)))\cdot\frac{\partial^2 F_\varepsilon}{\partial x^2}(t,X(t))\cdot g(t,X(t))^2,\end{eqnarray*}
$$\eta_\varepsilon\left[\frac{\partial g}{\partial t}\frac{\partial\beta_\varepsilon}{\partial\varepsilon}+\frac{\partial g}{\partial x}\cdot \phi+\frac{1}{2}\frac{\partial^2 g}{\partial x^2}\cdot\tilde{\phi}^2\right]+g\frac{\partial\eta_\varepsilon}{\partial\varepsilon}\qquad\qquad\qquad\qquad\qquad\qquad\qquad\qquad\qquad$$$$\qquad\qquad\qquad\qquad\qquad\qquad\qquad=\frac{\partial\phi}{\partial x}(t,F_\varepsilon(t,X(t)))\cdot\frac{\partial F_\varepsilon}{\partial x}(t,X(t))\cdot g(t,X(t))$$
\begin{eqnarray*}\eta_\varepsilon\frac{\partial g}{\partial x}\cdot\tilde{\phi}&=&\frac{\partial\tilde{\phi}}{\partial x}(t,F_\varepsilon(t,X(t)))\cdot\frac{\partial F_\varepsilon}{\partial x}(t,X(t))\cdot g(t,X(t)).\end{eqnarray*}

\vskip 0.5cm

Theese equations, when $\varepsilon=0$, provide our main theorem for this section.

\begin{teo}\label{teo1}
Let $\textbf{\emph{v}}=\textbf{\emph{v}}^D+\textbf{\emph{v}}^S$ be a vector field such that $\textbf{\emph{v}}^D=\uptau(t)\frac{\partial}{\partial t}+\phi(t,x)\frac{\partial}{\partial x}$ and $\textbf{\emph{v}}^S=\tilde{\phi}(t,x)\frac{\partial}{\partial x}$. If $\textbf{\emph{v}}$ is a stochastic symmetry for the SDE $$dX(t)=f(t,X(t))dt+g(t,X(t))dW(t),$$ then:
\begin{equation}\label{teo11}\begin{array}{rcl}\displaystyle\frac{\partial f}{\partial t}\uptau+\displaystyle\frac{\partial\uptau}{\partial t}f+\displaystyle\frac{\partial f}{\partial x}\cdot\phi+\displaystyle\frac{1}{2}\displaystyle\frac{\partial^2f}{\partial x^2}\cdot\tilde{\phi}^2&=&\displaystyle\frac{\partial \phi}{\partial t}+\displaystyle\frac{\partial \phi}{\partial x}\cdot f+ \displaystyle\frac{1}{2}\displaystyle\frac{\partial^2\phi}{\partial x^2}\cdot g^2\\[0.3cm]\displaystyle\frac{\partial f}{\partial x}\cdot\tilde{\phi}&=&\displaystyle\frac{\partial \tilde{\phi}}{\partial t}+\displaystyle\frac{\partial \tilde{\phi}}{\partial x}\cdot f+\displaystyle\frac{1}{2}\displaystyle\frac{\partial^2\tilde{\phi}}{\partial x^2}\cdot g^2\\[0.3cm]\displaystyle\frac{\partial g}{\partial t}\uptau+\displaystyle\frac{1}{2}\displaystyle\frac{\partial\uptau}{\partial t}g+\displaystyle\frac{\partial g}{\partial x}\cdot\phi+\displaystyle\frac{1}{2}\displaystyle\frac{\partial^2g}{\partial x^2}\cdot\tilde{\phi}^2&=&\displaystyle\frac{\partial \phi}{\partial x}\cdot g\\[0.3cm]\displaystyle\frac{\partial g}{\partial x}\cdot\tilde{\phi}&=&\displaystyle\frac{\partial \tilde{\phi}}{\partial x}\cdot g,\end{array}\end{equation}
where all functions are being evaluated at $(t,X(t))$ (except for $\uptau$ and its derivative, which are evaluated at $t$).
\end{teo}

\begin{obs}Note that
\begin{enumerate}
\item If we take $\tilde{\phi}=0$ in (\ref{teo11}), we obtain determinig equations (\ref{condicaosimetria1}) and (\ref{condicaosimetria2}).
\item If we take $g=0$ in (\ref{teo11}), we obtain determinig equations
\begin{equation}\begin{array}{rcl}\displaystyle\frac{\partial f}{\partial t}\uptau+\displaystyle\frac{\partial\uptau}{\partial t}f+\displaystyle\frac{\partial f}{\partial x}\cdot\phi+\displaystyle\frac{1}{2}\displaystyle\frac{\partial^2f}{\partial x^2}\cdot\tilde{\phi}^2&=&\displaystyle\frac{\partial \phi}{\partial t}+\displaystyle\frac{\partial \phi}{\partial x}\cdot f\\[0.3cm]\displaystyle\frac{\partial f}{\partial x}\cdot\tilde{\phi}&=&\displaystyle\frac{\partial \tilde{\phi}}{\partial t}+\displaystyle\frac{\partial \tilde{\phi}}{\partial x}\cdot f,\end{array}\end{equation}
which provides necessary conditions to obtain stochastic symmetries of deterministic ordinary differential equations.
\end{enumerate}\end{obs}

\begin{ex}
\emph{A simple example, is the equation of a Brownian motion in $\bbr$:} \begin{equation}\label{ex11}dX(t)=dW(t).\end{equation}
\end{ex}
Here, determining equations (\ref{condicaosimetria1}) and (\ref{condicaosimetria2}) are
\begin{eqnarray*}\label{}0&=&\phi_t+\frac{1}{2}\phi_{xx}\\\label{}\frac{1}{2}\uptau_t&=&\phi_x,\end{eqnarray*}
whose solution is 
\begin{eqnarray*}\label{}\uptau(t)&=&2c_1t+c_3\\\label{}\phi(t,x)&=&c_1x+c_2,\end{eqnarray*}
which generates the symmetries
\begin{equation}\label{ex12} X_1=2t\frac{\partial}{\partial t}+x\frac{\partial}{\partial x},\qquad X_2=\frac{\partial}{\partial x},\qquad X_3=\frac{\partial}{\partial t}.\end{equation}

For stochastic symmetries, determining equations (\ref{teo11}) are
$$\begin{array}{rcl}0&=&\displaystyle\frac{\partial \phi}{\partial t}+ \displaystyle\frac{1}{2}\displaystyle\frac{\partial^2\phi}{\partial x^2}\\[0.3cm]0&=&\displaystyle\frac{\partial \tilde{\phi}}{\partial t}+\displaystyle\frac{1}{2}\displaystyle\frac{\partial^2\tilde{\phi}}{\partial x^2}\\[0.3cm]\displaystyle\frac{1}{2}\displaystyle\frac{\partial\uptau}{\partial t}&=&\displaystyle\frac{\partial \phi}{\partial x}\\[0.3cm]0&=&\displaystyle\frac{\partial \tilde{\phi}}{\partial x},\end{array}$$
whose solution is
\begin{eqnarray*}\label{}\uptau(t)&=&2c_1t+c_3\\\label{}\phi(t,x)&=&c_1x+c_2\\\tilde{\phi}(t,x)&=&c_4.\end{eqnarray*}
These generates the same vector fields as in (\ref{ex12})
\begin{equation} X_1=\left[2t\frac{\partial}{\partial t}+x\frac{\partial}{\partial x}\right]^D,\qquad X_2=\left[\frac{\partial}{\partial x}\right]^D,\qquad X_3=\left[\frac{\partial}{\partial t}\right]^D\end{equation}
plus a new one
\begin{equation} X_4=\left[\frac{\partial}{\partial x}\right]^S.\end{equation}
So, in the case of the Brownian motion, stochastic symmetries produce a new vector field of symetry, that is $X_4$, which is a vector field with only stochastic part.

\begin{ex}
\emph{Consider the Langevin equation:} \begin{equation}\label{ex21}dX(t)=aX(t)dt+bdW(t).\end{equation}
\end{ex}
Determining equations (\ref{condicaosimetria1}) and (\ref{condicaosimetria2}) provide the symmetries
\begin{equation}\label{ex22} X_1=e^{at}\frac{\partial}{\partial x},\qquad X_2=\frac{e^{2at}}{a}\frac{\partial}{\partial t}+e^{2at}x\frac{\partial}{\partial x},\qquad X_3=\frac{\partial}{\partial t}.\end{equation}
On the other hand, determining equations (\ref{teo11}) provide the stochastic symmetries
\begin{equation} X_1=\left[e^{at}\frac{\partial}{\partial x}\right]^D,\qquad X_2=\left[\frac{e^{2at}}{a}\frac{\partial}{\partial t}+e^{2at}x\frac{\partial}{\partial x}\right]^D,\qquad X_3=\left[\frac{\partial}{\partial t}\right]^D\end{equation}
plus a new one
\begin{equation} X_4=\left[e^{at}\frac{\partial}{\partial x}\right]^S.\end{equation}
Again, stochastic symmetries produce a new vector field of symmetry with only stochastic part. However, this is not always the case, as the next example will show.

\begin{ex}
\emph{Consider the following equation in $\bbr$:} \begin{equation}\label{ex31}dX(t)=\frac{a}{X(t)}dt+dW(t).\end{equation}
\end{ex}
Symmetries for (\ref{ex31}) are
\begin{equation}\label{ex32} X_1=2t\frac{\partial}{\partial t}+x\frac{\partial}{\partial x},\qquad X_2=\frac{\partial}{\partial t}.\end{equation}
And the stochastic symmetries are just the same, that is,
\begin{equation} X_1=\left[2t\frac{\partial}{\partial t}+x\frac{\partial}{\partial x}\right]^D,\qquad X_2=\left[\frac{\partial}{\partial t}\right]^D.\end{equation}

\section{Transformations of SDEs}

Consider the SDE \begin{equation}\label{ede3}dX(t)=f(t,X(t))dt+g(t,X(t))dW(t),\end{equation} as before.  A possible way to obtain an explicit solution to (\ref{ede3}) is to transform it to a previously known differential equation. Suppose we are able to solve \begin{equation}\label{ede4}dY(s)=h(s,Y(s))ds+\sigma(s,Y(s))d\tilde{W}(s),\end{equation} where $h$ and $\sigma$ belongs to the same spaces as $f$ and $g$, respectively, and $\tilde{W}(s)$ is the Wiener process $W(t)$ after the change $t\mapsto s$. If one can find a function $\mu:(t,x)\mapsto(s,y)$ that transforms the SDE (\ref{ede3}) into the target SDE (\ref{ede4}), then one could obtain the solution to (\ref{ede3}) via the inverse of $\mu$.

What we are going to do is to provide a way to find such transformation $\mu$. In \cite{bluman}, G. W. Bluman et al. developed this theory for deterministic differential equations. We are going to extend it for stochastic differential equations.

Let $\textbf{v}$ be a vector field in the $(t,x)-$space such that its flow $(t,x)\mapsto(\overline{t},\overline{x})$ is a symetry for (\ref{ede3}) and let $\textbf{u}$ be vector field in the $(s,y)-$space whose flow $(s,y)\mapsto(s^*,y^*)$ is a symmetry for (\ref{ede4}). If $\mu(t,x)=(\mu_1(t,x),\mu_2(t,x))$ transforms (\ref{ede3}) in (\ref{ede4}), then the following diagram must comutate
$$\begin{array}{ccc}(t,x)&\longrightarrow&(\overline{t},\overline{x})\\\downarrow&&\downarrow\\(s,y)&\longrightarrow&(s^*,y^*).\end{array}$$
This means that we must have, for all $(t,x)$,
\begin{equation}\label{mu1}(\mu_1(\overline{t},\overline{x}),\mu_2(\overline{t},\overline{x}))=(\mu_1(t,x)^*,\mu_2(t,x)^*).\end{equation}

Now, let us rewrite equation (\ref{mu1}) using the coordinates of the flows associated to $\textbf{v}$ and $\textbf{u}$. Suppose $\textbf{v}=\textbf{v}^D+\textbf{v}^E$, with $\textbf{v}^D=\tau(t,x)\frac{\partial}{\partial t}+\phi(t,x)\frac{\partial}{\partial x}$ and $\textbf{v}^E=\tilde{\phi}(t,x)\frac{\partial}{\partial x}$ such that $(\overline{t},\overline{x})=(\beta_\varepsilon(t,x),F_\varepsilon(t,x))$. Suppose $\textbf{u}=\textbf{u}^D+\textbf{u}^E$, with $\textbf{u}^D=\rho(s,y)\frac{\partial}{\partial s}+\psi(s,y)\frac{\partial}{\partial y}$ and $\textbf{u}^E=\tilde{\psi}(s,y)\frac{\partial}{\partial y}$ such that $(s^*,y^*)=(b_\varepsilon(s,y),G_\varepsilon(s,y))$. Using this notatiton we rewrite equation (\ref{mu1}) as
\begin{equation}\label{mu2}(\mu_1(\beta_\varepsilon(t,x),F_\varepsilon(t,x)),\mu_2(\beta_\varepsilon(t,x),F_\varepsilon(t,x)))=(b_\varepsilon(\mu(t,x)),G_\varepsilon(\mu(t,x))),\end{equation}
from which we obtain 
\begin{equation}\label{mu3}\begin{array}{rcl}b_\varepsilon(\mu(t,x))&=&\mu_1(\beta_\varepsilon(t,x),F_\varepsilon(t,x))\\[0.3cm] G_\varepsilon(\mu(t,x))&=&\mu_2(\beta_\varepsilon(t,x),F_\varepsilon(t,x)).\end{array}\end{equation}
Differentiating equations (\ref{mu3}) with respect to $\varepsilon$ and comparing stochastic and deterministic terms, we get
\begin{equation}\label{mu4}\begin{array}{rcl}\displaystyle\frac{\partial b_\varepsilon}{\partial\varepsilon}\circ\mu&=&\displaystyle\frac{\partial\mu_1}{\partial t}\displaystyle\frac{\partial\beta_\varepsilon}{\partial\varepsilon}+\displaystyle\frac{\partial\mu_1}{\partial x}\phi+\displaystyle\frac{1}{2}\displaystyle\frac{\partial^2\mu_1}{\partial x^2}\tilde{\phi}^2\\[0.3cm] 0&=&\displaystyle\frac{\partial\mu_1}{\partial x}\tilde{\phi}\\[0.3cm] \psi\circ\mu&=&\displaystyle\frac{\partial\mu_2}{\partial t}\displaystyle\frac{\partial\beta_\varepsilon}{\partial\varepsilon}+\displaystyle\frac{\partial\mu_2}{\partial x}\phi+\displaystyle\frac{1}{2}\displaystyle\frac{\partial^2\mu_2}{\partial x^2}\tilde{\phi}^2\\[0.3cm] \tilde{\psi}\circ\mu&=&\displaystyle\frac{\partial\mu_2}{\partial x}\tilde{\phi}.\end{array}\end{equation}
Making $\varepsilon=0$ on equations (\ref{mu4}), we obtain our result.

\begin{teo}
If $\mu$ transforms a SDE \begin{equation}dX(t)=f(t,X(t))dt+g(t,X(t))dW(t)\end{equation} into the SDE \begin{equation}dY(s)=h(s,Y(s))ds+\sigma(s,Y(s))d\tilde{W}(s),\end{equation} then \begin{equation}\label{mu5}\begin{array}{rcl}\rho\circ\mu&=&\displaystyle\frac{\partial\mu_1}{\partial t}\tau+\displaystyle\frac{\partial\mu_1}{\partial x}\phi+\displaystyle\frac{1}{2}\displaystyle\frac{\partial^2\mu_1}{\partial x^2}\tilde{\phi}^2\\[0.3cm] 0&=&\displaystyle\frac{\partial\mu_1}{\partial x}\tilde{\phi}\\[0.3cm] \psi\circ\mu&=&\displaystyle\frac{\partial\mu_2}{\partial t}\tau+\displaystyle\frac{\partial\mu_2}{\partial x}\phi+\displaystyle\frac{1}{2}\displaystyle\frac{\partial^2\mu_2}{\partial x^2}\tilde{\phi}^2\\[0.3cm] \tilde{\psi}\circ\mu&=&\displaystyle\frac{\partial\mu_2}{\partial x}\tilde{\phi}\,.\end{array}\end{equation}
\end{teo}

\vspace{0.5cm}
Let us give an example of how to use this tool. In \cite{kozlov}, R. Kozlov gives a classification of ordinary SDEs in $\bbr$ in accordance to its symmetries. The way he does that is by using transformations between the SDEs. However, he does not mention how those transformations were obtained. Let us do so, using what we introduced in this section.

\begin{ex}
\begin{equation}\label{exmu1}\mu(t,x)=\left(-\frac{e^{-2\alpha t}}{2\alpha},e^{-\alpha t}\left(x+\frac{\beta}{\alpha}\right)\right)\end{equation}
\emph{transforms} \begin{equation}\label{ex11}dX(t)=\left(\alpha X(t)+\beta\right)dt+dW(t), \qquad \alpha\neq0\end{equation} \emph{into} \begin{equation}\label{ex12}dX(t)=dW(t).\end{equation}
\end{ex}

Symmetries of (\ref{ex11}) are generated by $$X_1=\frac{\partial}{\partial t},\qquad X_2=e^{2\alpha t}\frac{\partial}{\partial t}+(\alpha x+\beta)e^{2at}\frac{\partial}{\partial x},\qquad X_3=e^{\alpha t}\frac{\partial}{\partial x},$$ and symmetries of (\ref{ex12}) by $$Y_1=\frac{\partial}{\partial t},\qquad Y_2=2t\frac{\partial}{\partial t}+x\frac{\partial}{\partial x},\qquad Y_3=\frac{\partial}{\partial x}.$$
Note that their commutators are $$[X_1,X_2]=2\alpha X_2,\qquad[X_1,X_3]=\alpha X_3,\qquad[X_2,X_3]=0$$ and $$[Y_1,Y_2]=2Y_1,\qquad[Y_1,Y_3]=0,\qquad[Y_2,Y_3]=Y_3.$$ First, we adjust the coefficients. Replacing $X_i$ for $\tilde{X}_i=a_1^iX_1+a_2^iX_2+a_3^iX_3$ and requiring that $$[\tilde{X}_1,\tilde{X}_2]=2\tilde{X}_1,\qquad[\tilde{X}_1,\tilde{X}_3]=0,\qquad[\tilde{X}_2,\tilde{X}_3]=\tilde{X}_3,$$ we obtain $a^1_2=-1, a^2_1=\frac{1}{\alpha},a^3_3=1$ and $a^i_j=0$ in the other cases, which means that $$\tilde{X}_1=-X_2,\qquad\tilde{X}_2=\frac{1}{\alpha}X_1,\qquad\tilde{X}_3=X_3.$$ Hence, $\tilde{X}_1,\tilde{X}_2,\tilde{X}_3$ generates the same Lie algebra as $X_1,X_2,X_3$ and their commutator agree with the commutators of $Y_1,Y_2,Y_3$. Then, to find a transformation $\mu$ that takes (\ref{ex11}) into (\ref{ex12}), conditions (\ref{mu5}) must be satisfied for each pair of symmetries $\tilde{X}_i,Y_i$. This implies that $\mu$ must satisfy the following system: 
\begin{equation}\label{exsys1}\left\{\begin{array}{rcl}1&=&-e^{2\alpha t}\frac{\partial\mu_1}{\partial t}+(\alpha x+\beta)e^{2\alpha t}\frac{\partial\mu_1}{\partial x}\\[0.3cm]0&=&-e^{2\alpha t}\frac{\partial\mu_2}{\partial t}+(\alpha x+\beta)e^{2\alpha t}\frac{\partial\mu_2}{\partial x}\\[0.3cm]2\mu_1&=&\frac{1}{\alpha}\frac{\partial\mu_1}{\partial t}\\[0.3cm]\mu_2&=&\frac{1}{\alpha}\frac{\partial\mu_2}{\partial t}\\[0.3cm]0&=&e^{\alpha t}\frac{\partial\mu_1}{\partial x}\\[0.3cm]1&=&e^{\alpha t}\frac{\partial\mu_2}{\partial x}.\end{array}\right.\end{equation} Solving system (\ref{exsys1}), one finds exactly (\ref{exmu1}).

\vspace{0.5cm}
Note that the other transformations given in \cite{kozlov} can be obtained in the same way.



\begin{thebibliography}{99}

\bibitem{albeverio} Albeverio, S. and Fei, S-M.: \emph{A remark on symmetry of stochastic dynamical systems and their conserved quantities.} J. Phys. A \textbf{28}, no. 22, 6363--6371, 1995.

\bibitem{bluman1} Bluman, G. W. and Anco, S. C.: \emph{Symmetry and integration methods for differential equations.} Applied Mathematical Sciences \textbf{154}, Springer-Verlag, New York, 2002.

\bibitem{bluman} Bluman, G. W., Cheviakov, A. F. and Anco, S. C.: \emph{Applications of symmetry methods to partial differential equations.} Applied Mathematical Sciences, vol. 168, Springer, New York, 2010.

\bibitem{gaeta} Gaeta, G. and Quintero, N. R.: \emph{Lie-point symmetries and differential equations.} J. Phys. A: Math. Gen. \textbf{32}, 8425--8505, 1999.

\bibitem{kozlov} Kozlov, R.: \emph{The group classification of a scalar stochastic differential equation}. J. Phys. A: Math. Theor. \textbf{43}, no. 5, 2010.

\bibitem{olver} Olver, P. J.: \emph{Applications of Lie groups to differential equations.} Graduate Texts in Mathematics, vol. 107, Springer-Verlag, New York, 1993.

\bibitem{oksendal} Oksendal, B.: \emph{Stochastic differential equations. An introduction with applications.} Fifth edition, Universitext, Springer-Verlag, Berlin, 1998.

\bibitem{meleshko} Srihirun, B., Meleshko, S. and Schulz, E.: \emph{On the definition of an admitted Lie group for stochastic differential equations.} Commun. Nonlinear Sci. Numer. Simul. \textbf{12}, 1379--1389, 2007.

\bibitem{stefani} Stephani, H.: \emph{Differential equations: their solution using symmetries.} Cambridge University Press, 1989.
\bibitem{unal} Ünal, G.: \emph{Symmetries of Itô and Stratonovich dynamical systems and their conserved quantities.} Nonlinear Dyn. \textbf{32}, 417--426, 2003.

\bibitem{wafo} Wafo Soh, C. and Mahomed, F. M.: \emph{Integration of stohastic ordinary differential equations from a symmetry standpoint.} J. Phys. A: Math. Gen. \textbf{34}, 177--192, 2001.


\end{thebibliography}
\end{document}